\documentclass[11pt]{article}
\usepackage{amsfonts}
\usepackage{amsmath}
\usepackage{amssymb}
\usepackage{graphicx}
\usepackage[dvipsnames,usenames]{color}
\usepackage[pdfstartview=FitH ]{hyperref}
\usepackage[hyperpageref]{backref}
\begin{document}

\baselineskip 18pt
\def\o{\over}
\def\e{\varepsilon}
\title{\Large\bf  A\ note\ on\ the\ number\ of\ coefficients\\
of\ automorphic\ $L-$functions\\
for\ $GL_m$\ with\ same\ signs}
\author{Chaohua\ \ Jia}
\date{}
\maketitle {\small \noindent {\bf Abstract.} Let $\pi$ be an
irreducible unitary cuspidal representation of $GL_m({\Bbb A}_{\Bbb
Q})$ and $L(s,\,\pi)$ be the global $L-$function attached to $\pi$.
If ${\rm Re}(s)>1$, $L(s,\,\pi)$ has a Dirichlet series expression.
When $\pi$ is self-contragradient, all the coefficients of Dirichlet
series are real. In this note, we shall give non-trivial lower
bounds for the number of positive and negative coefficients
respectively, which is an improvement on the recent work of Jianya
Liu and Jie Wu.

}

\vskip.3in
\noindent{\bf 1. Introduction}

Let $m\geq 2$ be an integer, $\pi=\otimes\pi_p$ be an irreducible
unitary cuspidal representation of $GL_m({\Bbb A}_{\Bbb Q})$. For
${\rm Re}(s)>1$, the global $L-$function is defined as
$$
L(s,\,\pi):=\prod_{p<\infty}L_p(s,\,\pi_p), \eqno (1.1)
$$
where
$$
L_p(s,\,\pi_p):=\prod_{1\leq j\leq m}\Bigl(1-{\alpha_\pi(p,\,j)\o
p^s}\Bigr)^{-1}. \eqno (1.2)
$$

To estimate the local parameters $\{\alpha_\pi(p,\,j)\}_{j=1}^m$
associated with $\pi_p$ is important for the study of automorphic
$L-$functions. The best known results are that for all primes $p$
and $1\leq j\leq m$,
$$
\alpha_\pi(p,\,j)\ll p^{\theta_m},     \eqno (1.3)
$$
where
$$
\theta_2={7\o 64},\qquad \theta_3={5\o 14},\qquad \theta_4={9\o 22}
\eqno (1.4)
$$
which are due to H. H. Kim and P. Sarnak[3],
$$
\theta_m={1\o 2}-{2\o m^2+1},\qquad\quad m\geq 5 \eqno (1.5)
$$
which are due to Wenzhi Luo, Z. Rudnick and P. Sarnak[6].

H. Jacquet and J. A. Shalika[2] proved that for ${\rm Re}(s)>1$, the
Euler product for $L(s,\,\pi)$ in (1.1) converges absolutely so that
we can write
$$
L(s,\,\pi)=\sum_{n=1}^\infty{\lambda_\pi(n)\o n^s}, \eqno (1.6)
$$
where
$$
\lambda_\pi(n)=\prod_{p^\nu\|\,n}\Bigl(\sum_{\nu_1+\cdots+\nu_m=
\nu}\alpha_\pi(p,\,1)^{\nu_1}\cdots\alpha_\pi(p,\,m)^{\nu_m}\Bigr).
\eqno (1.7)
$$

When $\pi$ is self-contragradient, $\lambda_\pi(n)$ is real for all
$n\geq 1$ so that one would be interested in the problem of sign
changes of $\lambda_\pi(n)$. Yan Qu[8] proved that, if $\pi$ is a
self-contragradient irreducible unitary cuspidal representation of
$GL_m({\Bbb A}_{\Bbb Q})$, then there must be infinitely many sign
changes in the sequence $\{\lambda_\pi(n)\}_{n=1}^\infty$, i.e.,
there are infinitely many $n$ such that $\lambda_\pi(n)>0$, and
there are infinitely many $n$ such that $\lambda_\pi(n)<0$.

Recently, Jianya Liu and Jie Wu[4] gave a quantitative version of
the above result of Yan Qu[8]. Write
$$
{\cal N}_\pi^+(x):=\sum_{\substack{n\leq x\\ \lambda_\pi(n)>0}}1,
\eqno (1.8)
$$
and
$$
{\cal N}_\pi^-(x):=\sum_{\substack{n\leq x\\ \lambda_\pi(n)<0}}1.
\eqno (1.9)
$$
They[4] showed that, if $\pi$ is a self-contragradient irreducible
unitary cuspidal representation of $GL_m({\Bbb A}_{\Bbb Q})$ and
$\theta_m$ is in (1.4) or (1.5), then for $x\geq x_0(\pi)$, one has
$$
{\cal N}_\pi^\pm(x)\gg_\pi x^{1-2\theta_m}(\log x)^{{2\o m}-m-4},
\eqno(1.10)
$$
unconditionally for $2\leq m\leq 4$ and under the Hypothesis H for
$m\geq 5$. The well known Hypothesis H is due to Z. Rudnick and P.
Sarnak[10], which is stated as follows.

{\bf Hypothesis H}. Let
$$
a_\pi(n):=
\begin{cases}
\alpha_\pi(p,\,1)^\nu+\cdots+\alpha_\pi(p,\,m)^\nu,\quad &{\rm if}\
n=p^\nu,\\
\qquad\qquad\quad 0, &{\rm otherwise}
\end{cases}
\eqno (1.11)
$$
for all primes $p$ and integers $\nu\geq 1$. Then for any fixed
integer $\mu\geq 2$,
$$
\sum_p{|a_\pi(p^\mu)|^2\log^2p\o p^\mu}<\infty. \eqno (1.12)
$$

In this note, we shall give another proof of (1.10) with an
improvement on the logarithm factor.

{\bf Theorem}. Let $\pi$ be a self-contragradient irreducible
unitary cuspidal representation of $GL_m({\Bbb A}_{\Bbb Q})$ and
$\theta_m$ be in (1.4) or (1.5). Then for $x\geq x_0(\pi)$, we have
$$
{\cal N}_\pi^\pm(x)\gg_\pi x^{1-2\theta_m}(\log x)^{-2[{m\theta_m\o
2}]-4}, \eqno(1.13)
$$
unconditionally for $2\leq m\leq 4$ and under the Hypothesis H for
$m\geq 5$.

It is easy to verify
$$
-2\Bigl[{m\theta_m\o 2}\Bigr]-4>{2\o m}-m-4,
$$
so that we can get an improvement on the logarithm factor for
(1.10).

\vskip.3in \noindent{\bf 2. Some lemmas}

{\bf Lemma 1}. Let $L(s,\,\pi)$ be defined as in (1.1) and (1.2).
Then $L(s,\,\pi)$ is an entire function. For real numbers
$\sigma\geq 0$ and $t$, we have
$$
L(\sigma+it,\,\pi)\ll_{\pi,\,\e}(|t|+1)^{\max({m\o
2}(1-\sigma),\,0)+\e}, \eqno (2.1)
$$
where $\e$ is a sufficiently small positive constant.

Proof. By the exposition in page 85 of [7], we know that
$L(s,\,\pi)$ is an entire function. The application of convexity
bound for $L(s,\,\pi)$ of G. Harcos[1] yields the estimate (2.1).

{\bf Lemma 2}. Let $\theta_m$ be in (1.4) or (1.5). We have
$$
\lambda_\pi(n)\ll_\pi n^{\theta_m+\e}, \eqno (2.2)
$$
and
$$
\lambda_\pi(p)=\alpha_\pi(p,\,1)+\cdots+\alpha_\pi(p,\,m)\ll_\pi
p^{\theta_m}. \eqno (2.3)
$$

They come from (1.3) and (1.7).

{\bf Lemma 3}. We have
$$
\sum_{n\leq x}\lambda_\pi(n)\Bigl(\log{x\o n}\Bigr)^{[{m\theta_m\o
2}]+1}\ll_{\pi,\,\e}x^{1-\theta_m-\e}, \eqno (2.4)
$$
where $\e$ is a sufficiently small positive constant.

Proof. For any positive integer $k$,
$$
{1\o 2\pi i}\int_{2-i\infty}^{2+i\infty}{y^s\o s^{k+1}}ds=
\begin{cases} {1\o k!}\log^ky,\qquad &{\rm if}\quad y\geq 1,\\
\quad 0,&{\rm if}\ 0<y<1.
\end{cases}
$$
Thus
\begin{align*}
&\ \,\sum_{n\leq x}\lambda_\pi(n)\Bigl(\log{x\o
n}\Bigr)^{[{m\theta_m\o 2}]+1}\\
&={1\o 2\pi i}\Bigl(\Bigl[{m\theta_m\o
2}\Bigr]+1\Bigr)!\int_{2-i\infty}^{2+i\infty}L(s,\,\pi){x^s\o
s^{[{m\theta_m\o 2}]+2}}ds.
\end{align*}

We move the line of integration to ${\rm Re}(s)=1-\theta_m-\e$. An
application of Lemma 1 produces
\begin{align*}
&\ \,\int_{2-i\infty}^{2+i\infty}L(s,\,\pi){x^s\o s^{[{m\theta_m\o
2}]+2}}ds\\
&=\int_{1-\theta_m-\e-i\infty}^{1-\theta_m-\e+i\infty}L(s,\,\pi){x^s\o
s^{[{m\theta_m\o 2}]+2}}ds\\
&\ll\int_{-\infty}^\infty (|t|+1)^{{m\o
2}(\theta_m+\e)+\e}{x^{1-\theta_m-\e}\o (|t|+1)^{[{m\theta_m\o
2}]+2}}dt\\
&=x^{1-\theta_m-\e}\int_{-\infty}^\infty{dt\o
(|t|+1)^{2-\{{m\theta_m\o 2}\}-{m\e\o 2}-\e}}\\
&\ll x^{1-\theta_m-\e}.
\end{align*}
Therefore
$$
\sum_{n\leq x}\lambda_\pi(n)\Bigl(\log{x\o n}\Bigr)^{[{m\theta_m\o
2}]+1}\ll_{\pi,\,\e}x^{1-\theta_m-\e}.
$$

So far the proof of Lemma 3 is finished.

{\bf Lemma 4}. We have
$$
\sum_{n\leq x}|\lambda_\pi(n)|^2\ll_\pi x.  \eqno (2.5)
$$

Proof. When $m=2$, R. A. Rankin[9] showed that
$$
\sum_{n\leq x}|\lambda_\pi(n)|^2=cx+O_\pi(x^{3\o 5})\ll_\pi x,
$$
where $c$ is a positive constant.

When $m\geq 3$, let $L(s,\,\pi\times\tilde{\pi})$ be the
Rankin-Selberg $L-$function associated to $\pi$ and its
contragradient $\tilde{\pi}$, which is defined as
$$
L(s,\,\pi\times\tilde{\pi}):=\prod_{p<\infty}L(s,\,\pi_p\times\tilde{\pi}_p).
$$
Write
$$
L(s,\,\pi\times\tilde{\pi})=\sum_{n=1}^\infty{\lambda_{\pi\times\tilde{\pi}}(n)\o
n^s}.
$$

By (3.10) in page 2884 of [5], we have
$$
\sum_{n\leq x}\lambda_{\pi\times\tilde{\pi}}(n)=c_\pi
x+O_{\pi,\,\e}(x^{{m^2-1\o m^2+1}+\e})\ll_\pi x,
$$
where $c_\pi$ is a positive constant. The discussion in page 2885 of
[5] yields
$$
|\lambda_\pi(n)|^2\leq\lambda_{\pi\times\tilde{\pi}}(n).
$$
Hence,
$$
\sum_{n\leq x}|\lambda_\pi(n)|^2\leq\sum_{n\leq
x}\lambda_{\pi\times\tilde{\pi}}(n)\ll_\pi x.
$$

So far the proof of Lemma 4 is finished.

{\bf Lemma 5}. Let $\pi$ be a self-contragradient irreducible
unitary cuspidal representation of $GL_m({\Bbb A}_{\Bbb Q})$. Then
there is a positive constant $c=c(\pi)$ such that
$$
\sum_{p\leq x}|\lambda_\pi(p)|^2\log p=x+O_\pi(xe^{-c\sqrt{\log x}})
\eqno (2.6)
$$
holds true unconditionally for $2\leq m\leq 4$ and under the
Hypothesis H for $m\geq 5$.

This is Theorem 3 of [11].

\vskip.3in
\noindent{\bf 3. The proof of Theorem}

We only prove for ${\cal N}_\pi^+(x)$. The proof for  ${\cal
N}_\pi^-(x)$ is same.

For $x\geq x_0(\pi)$, by (2.3) and Lemma 5, we have
\begin{align*}
&\ \,\sum_{n\leq x}|\lambda_\pi(n)|\Bigl(\log{x\o
n}\Bigr)^{[{m\theta_m\o 2}]+1}\\
&\geq\sum_{n\leq {x\o 2}}|\lambda_\pi(n)|\Bigl(\log{x\o
n}\Bigr)^{[{m\theta_m\o 2}]+1}\\
&\geq(\log 2)^{[{m\theta_m\o 2}]+1}\sum_{n\leq {x\o
2}}|\lambda_\pi(n)|\\
&\gg\sum_{p\leq {x\o 2}}|\lambda_\pi(p)|\\
&\gg{1\o x^{\theta_m}\log x}\sum_{p\leq {x\o 2}}|\lambda_\pi(p)|^2
\log p\\
&\gg{x^{1-\theta_m}\o \log x}.
\end{align*}
By Lemma 3,
$$
\sum_{n\leq x}\lambda_\pi(n)\Bigl(\log{x\o n}\Bigr)^{[{m\theta_m\o
2}]+1}\ll x^{1-\theta_m-\e}.
$$
Hence,
$$
\sum_{n\leq x}{|\lambda_\pi(n)|+\lambda_\pi(n)\o 2}\Bigl(\log{x\o
n}\Bigr)^{[{m\theta_m\o 2}]+1}\gg{x^{1-\theta_m}\o \log x}. \eqno
(3.1)
$$

On the other hand, we note that $\lambda_\pi(n)$ is real for all
$n\geq 1$ and use the Cauchy's inequality and Lemma 4 to get
\begin{align*}
&\ \,\sum_{n\leq x}{|\lambda_\pi(n)|+\lambda_\pi(n)\o
2}\Bigl(\log{x\o n}\Bigr)^{[{m\theta_m\o 2}]+1}\\
&=\sum_{\substack{n\leq x\\ \lambda_\pi(n)>0}}|\lambda_\pi(n)|\Bigl(
\log{x\o n}\Bigr)^{[{m\theta_m\o 2}]+1}\\
&\leq(\log x)^{[{m\theta_m\o 2}]+1}\sum_{\substack{n\leq x\\
\lambda_\pi(n)>0}}|\lambda_\pi(n)|\\
&\leq(\log x)^{[{m\theta_m\o 2}]+1}\Bigl(\sum_{n\leq
x}|\lambda_\pi(n)|^2\Bigr)^{1\o 2}\Bigl(\sum_{\substack{n\leq x\\
\lambda_\pi(n)>0}}1\Bigr)^{1\o 2}\\
&\ll x^{1\o 2}(\log x)^{[{m\theta_m\o 2}]+1}({\cal N}_\pi^+(x))^{1\o
2}.
\end{align*}
Thus
\begin{align*}
{x^{1-\theta_m}\o \log x}&\ll x^{1\o 2}(\log x)^{[{m\theta_m\o
2}]+1}({\cal N}_\pi^+(x))^{1\o 2},\\
{\cal N}_\pi^+(x)&\gg x^{1-2\theta_m}(\log x)^{-2[{m\theta_m\o
2}]-4}.
\end{align*}
Therefore the proof of Theorem is complete.

\vskip.3in
\noindent{\bf Acknowledgements}

The author would like to thank his colleagues and friends in the
``ergodic prime number theorem 2014'' seminar in the Morningside
Center of Chinese Academy of Sciences, especially Dr. Hengcai Tang
and Dr. Qinghua Pi, for their helpful discussion.

The author is supported by the National Key Basic Research Program
of China (Project No. 2013CB834202) and the National Natural Science
Foundation of China (Grant No. 11371344, Grant No. 11321101).


\bigskip

\

\

Chaohua Jia

Institute of Mathematics, Academia Sinica, Beijing 100190, P. R.
China

Hua Loo-Keng Key Laboratory of Mathematics, Chinese Academy of
Sciences, Beijing 100190, P. R. China

E-mail: jiach@math.ac.cn


\vskip.6in

\begin{thebibliography}{9}

\bibitem{1} G. Harcos, {\it Uniform approximate functional equation
for principal $L-$functions}, Inter. Math. Research Notes, {\bf
13}(2002), 923-932; ibid. {\bf 18}(2004), 659-660.

\bibitem{2} H. Jacquet and J. A. Shalika, {\it On Euler products
and the classification of automorphic representations I}, Amer. J.
Math., {\bf 103}(1981), 499-558.

\bibitem{3} H. H. Kim and P. Sarnak, {\it Refined estimates towards
the Ramanujan and Selberg conjectures}, J. Amer. Math. Soc., {\bf
16}(2003), no.1, 175-183.

\bibitem{4} Jianya Liu and Jie Wu, {The number of coefficients of
automorphic $L-$functions for $GL_m$ of same signs},
arXiv:1404.6867v1[math.NT].

\bibitem{5} Guangshi L\"u, {\it On sums involving coefficients of
automorphic $L-$functions}, Proc. Amer. Math. Soc., {\bf 137}(2009),
2879-2887.

\bibitem{6} Wenzhi Luo, Z. Rudnick and P. Sarnak, {\it On the
generalized Ramanujan conjecture for $GL(n)$}, Proceedings of
Symposia in Pure Mathematics, vol.{\bf 66}, part 2, 1999, 301-310.

\bibitem{7} Yan Qu, {\it The prime number theorem for automorphic
$L-$functions for $GL_m$}, J. Number Theory, {\bf 122}(2007), 84-99.

\bibitem{8} Yan Qu, {\it Linnik-type problems for automorphic
$L-$functions}, J. Number Theory, {\bf 130}(2010), 786-802.

\bibitem{9} R. A. Rankin, {\it Contributions to the theory of
Ramanujan's function $\tau(n)$ and similar arithmetical functions
II. The order of the Fourier coefficients of the integral modular
forms}, Proc. Cambridge Phil. Soc., {\bf 35}(1939), 357-372.

\bibitem{10} Z. Rudnick and P. Sarnak, {\it Zeros of principal $L-$functions
and random matrix theory}, Duke Math. J., {\bf 81}(1996), 269-322.

\bibitem{11} Jie Wu and Yangbo Ye, {\it Hypothesis H and the prime
number theorem for automorphic representations}, Functiones et
Approximatio, {\bf 37}(2007), 461-471.

\end{thebibliography}
\end{document}